\documentclass[preprint,12pt]{elsarticle}

\usepackage{amssymb}
\usepackage{amsthm}
\usepackage{amsmath}
\usepackage{fancybox}

\newtheorem{thm}{Theorem}
\newtheorem{lem}{Lemma}
\newdefinition{rem}{Remark}
\newdefinition{exa}{Example}
\newproof{pf}{{\bf Proof}}
\newproof{pfirst}{{\bf Proof of Theorem~\protect\ref{t:first}}}
\newcommand{\eq}[1]{\mbox{\rm(\ref{#1})}}
\newcommand{\rra}{\rightrightarrows}

\catcode`\@=11  \@addtoreset{equation}{section}  \catcode`\@=12
\newcommand{\RB}{\mathbb{R}}   \newcommand{\NB}{\mathbb{N}}
\newcommand{\DC}{\mathcal{D}}
\newcommand{\vep}{\varepsilon}

\journal{arXiv (math.FA)\quad}

\begin{document}

\begin{frontmatter}

\title{The joint modulus of variation of metric space valued functions
       and pointwise selection principles} %\tnoteref{t1}}

\author[hsenn]{Vyacheslav V.~Chistyakov\corref{cor1}}
\ead{czeslaw@mail.ru, vchistyakov@hse.ru}

\author[hsenn]{Svetlana A.~Chistyakova}
\ead{schistyakova@hse.ru}

\cortext[cor1]{Corresponding author.}

\address[hsenn]{Department of Informatics, Mathematics and Computer Science,\\
National Research University Higher School of Economics,\\
Bol'shaya Pech{\"e}rskaya Street 25/12, Nizhny Novgorod\\
603155, Russian Federation}

\begin{abstract}
Given $T\subset\RB$ and a metric space $M$, we introduce a nondecreasing sequence
of pseudometrics $\{\nu_n\}$ on $M^T$ (the set of all functions from $T$ into $M$),
called the {joint modulus of variation}. We prove that {if two sequences
of functions $\{f_j\}$ and $\{g_j\}$ from $M^T$ are such that $\{f_j\}$ is pointwise
precompact, $\{g_j\}$~is pointwise convergent, and the limit superior of $\nu_n(f_j,g_j)$
as $j\to\infty$ is $o(n)$ as $n\to\infty$, then $\{f_j\}$ admits a pointwise convergent
subsequence whose limit is a conditionally regulated function}. We illustrate the sharpness
of this result by examples (in particular, the assumption on the $\limsup$ is necessary
for uniformly convergent sequences $\{f_j\}$ and $\{g_j\}$, and `almost necessary'
when they converge pointwise) and show that most of the known Helly-type pointwise
selection theorems are its particular cases.
\end{abstract}

\begin{keyword}
joint modulus of variation \sep metric space \sep regulated function \sep
pointwise convergence \sep selection principle \sep generalized variation

{\em MSC\,2000:} 26A45 \sep 28A20 \sep 54C35 \sep 54E50
\end{keyword}
\end{frontmatter}

\section{Introduction} \label{s:intro}

The purpose of this paper is to present a new sufficient condition (which is almost
necessary) on a pointwise precompact sequence $\{f_j\}\equiv\{f_j\}_{j=1}^\infty$
of functions $f_j$ mapping a subset $T$ of the real line $\RB$ into a metric space~$(M,d)$,
under which the sequence admits a pointwise convergent subsequence.
The historically first result in this direction is the classical \emph{Helly Selection Principle},
in which the assumptions are as follows: $T=[a,b]$ is a closed interval, $M=\RB$, and
$\{f_j\}$ is uniformly bounded and consists of monotone functions (\cite{Helly},
\cite[II.8.9--10]{Hild}, \cite[VIII.4.2]{Nat}, and \cite[Theorem~1.3]{Sovae} if
$T\subset\RB$ is arbitrary). Since a real function on $T$ of bounded (Jordan) variation
is the difference of two nondecreasing bounded functions, Helly's theorem extends
to uniformly bounded sequences of functions, whose Jordan's variations are uniformly
bounded. Further generalizations of the latter pointwise selection principle are concerned
with replacement of Jordan's variation by more general notions of variation
(\cite{jmaa00,Bhakta}, \cite{JDCS}--\cite{Sovae}, \cite{ChG1,ChG2,CyKe,Fr,Gnilka,
Ko1,Ko2,MuOr,Schramm,Wat76}). In all these papers, the pointwise limit of the extracted
subsequence of $\{f_j\}$ is a function of bounded generalized variation (in the
corresponding sense), and so, it is a regulated function (with finite one-sided limits
at all points of the domain). Note that pointwise selection principles (or the sequential
compactness in the topology of pointwise convergence) and regulated functions are
of importance in real analysis (\cite{GNW,Hild,Nat}), stochastic analysis and generalized
integration (\cite{Muld}), optimization (\cite{Barbu,Mordu}), set-valued analysis
(\cite{jmaa00,Sovae,JFA05,jmaa07,Hermes}), and other fields.

\smallbreak
A unified approach to the diverse selection principles mentioned above was proposed
in \cite{jmaa05,Ischia}. It is based on the notion of \emph{modulus of variation\/}
of a function introduced in \cite{Chan74,Chan75} (see also \cite[11.3.7]{GNW})
and does not refer to the uniform boundedness of variations of any kind, and so,
can be applied to sequences of non-regulated functions. However, the pointwise limit
of the extracted subsequence of $\{f_j\}$ is \emph{again\/} a regulated function.
In order to clarify this situation and expand the amount of sequences having pointwise
convergent subsequences, we define the notion of the \emph{joint
modulus of variation\/} for metric space valued functions: this is a certain sequence
of pseudometrics $\{\nu_n\}$ on the product set $M^T$ (of all functions from $T$
into $M$). Making use of $\{\nu_n\}$, we obtain a powerful pointwise selection principle
(see Theorem~\ref{t:first} in Section~\ref{s:main}). Putting $g_j=c$ for all $j\in\NB$,
where $c:T\to M$ is a constant function, we get the selection principle from \cite{jmaa05},
which already contains all selection principles alluded to above as particular cases.
In contrast to results from \cite{jmaa05,Ischia}, the pointwise limit $f$ from
Theorem~\ref{t:first} may not be regulated in general---this depends on the limit
function $g$, namely, since $\nu_n(f,g)=o(n)$, the function $f$ is only conditionally
regulated with respect to~$g$ (in short, \emph{$g$-re\-gu\-la\-ted\/}). In particular, if
$g=c$, \mbox{then $f$ is regulated in the usual sense.}

Finally, we point out that following the ideas from \cite{MatTr}, Theorem~\ref{t:first} may
be extended to sequences of functions with values in a \emph{uniform space\/}~$M$.

\smallbreak
The paper is organized as follows. In Section~\ref{s:main}, we present necessary
definitions and our main result, Theorem~\ref{t:first}. In Section~\ref{s:joint},
we establish essential properties of the joint modulus of variation, which are needed
in the proof of Theorem~\ref{t:first} in Section~\ref{s:proof}. Section~\ref{s:reg} is
devoted to the study of $g$-regulated (and, in particular, regulated) functions.
In the final Section~\ref{s:ext}, we extend the Helly-type selection theorems from
\cite{Fr} and \cite{CyKe,Ko1} by exploiting Theorem~\ref{t:first}.

\section{Main result} \label{s:main}
 
Let $\varnothing\ne T\subset\RB$, $(M,d)$ be a metric space with metric $d$, and
$M^T$ denote the set of all functions $f:T\to M$ mapping $T$ into $M$. The letter~$c$
stands for a \emph{constant\/} function $c:T\to M$ (i.e., $c(s)=c(t)$ in $M$
for all $s,t\in T$).

\smallbreak
The \emph{joint oscillation\/} of two functions $f,g\in M^T$ is the quantity
  $$|(f,g)(T)|=\sup\bigl\{|(f,g)(\{s,t\})|:s,t\in T\bigr\}\in[0,\infty],$$
where
  \begin{equation} \label{e:join}
|(f,g)(\{s,t\})|=\sup_{z\in M}\bigl|d(f(s),z)+d(g(t),z)-d(f(t),z)-d(g(s),z)\bigr|
  \end{equation}
is the \emph{joint increment\/} of $f$ and $g$ on the two-point set $\{s,t\}\subset T$,
for which the following two inequalities hold:
  \begin{align}
|(f,g)(\{s,t\})|&\le d(f(s),f(t))+d(g(s),g(t)),\label{e:fg1}\\[2pt]
|(f,g)(\{s,t\})|&\le d(f(s),g(s))+d(f(t),g(t)).\label{e:fg2}
  \end{align}

Since $|(f,c)(\{s,t\})|=d(f(s),f(t))$ ($=$\,the increment of $f$ on $\{s,t\}\subset T$) is
independent of $c$, the quantity $|f(T)|=|(f,c)(T)|$ is the usual \emph{oscillation\/}
of $f$ on $T$, also known as the \emph{diameter of the image\/}
$f(T)=\{f(t):t\in T\}\subset M$. Clearly, by \eq{e:fg1}, $|(f,g)(T)|\le|f(T)|+|g(T)|$.

\smallbreak
We denote by $\mbox{\rm B}(T;M)=\{f\in M^T:|f(T)|<\infty\}$ the family of all
\emph{bounded\/} functions on $T$ equipped with the \emph{uniform metric\/}
$d_\infty$ given by
  $$d_\infty(f,g)=\sup_{t\in T}d(f(t),g(t))\quad\,\mbox{for}\quad\,
     f,g\in\mbox{\rm B}(T;M)$$
($d_\infty$ is an \emph{extended\/} metric on $M^T$, i.e., may assume the value~%
$\infty$). We have
  $$d_\infty(f,g)\le d(f(s),g(s))+|f(T)|+|g(T)|\quad\mbox{for \,all}\quad s\in T$$
and, by virtue of \eq{e:fg2}, $|(f,g)(T)|\le2d_\infty(f,g)$.

\smallbreak
If $n\in\NB$, we write $\{I_i\}_1^n\prec T$ to denote a collection of $n$ two-point
subsets $I_i=\{s_i,t_i\}$ of $T$ ($i=1,\dots,n$) such that
$s_1<t_1\le s_2<t_2\le\dots\le s_{n-1}<t_{n-1}\le s_n<t_n$
(so that the intervals $[s_1,t_1],\dots,[s_n,t_n]$ with end-points in $T$ are
non-overlapping). We say that a collection $\{I_i\}_1^n\prec T$ with $I_i=\{s_i,t_i\}$ is a
\emph{partition\/} of $T$ if (setting $t_0=s_1$) $s_i=t_{i-1}$ for all $i=1,\dots,n$,
which is written as $\{t_i\}_0^n\prec T$.

\smallbreak
The \emph{joint modulus of variation\/} of two functions $f,g\in M^T$ is the sequence
$\{\nu_n(f,g)\}_{n=1}^\infty\subset[0,\infty]$ defined by
  \begin{equation} \label{e:jmv}
\nu_n(f,g)=\sup\,\biggl\{\sum_{i=1}^n|(f,g)(I_i)|:\{I_i\}_1^n\prec T\biggr\}\quad\,
\mbox{for \,all}\quad\,n\in\NB,
  \end{equation}
where $|(f,g)(I_i)|=|(f,g)(\{s_i,t_i\})|$ is the quantity from \eq{e:join} if $I_i=\{s_i,t_i\}$
(for finite $T$ with the number of elements $\#(T)\ge2$, we make use of \eq{e:jmv}
for $n\le\#(T)-1$, and set $\nu_n(f,g)=\nu_{\#(T)-1}(f,g)$ for all $n>\#(T)-1$).

\smallbreak
Note that, given $f,g\in M^T$, we have $\nu_1(f,g)=|(f,g)(T)|$ and
  \begin{equation} \label{e:nuon}
\nu_1(f,g)\le\nu_n(f,g)\le n\nu_1(f,g)\quad\,\mbox{for \,all}\quad\,n\in\NB.
  \end{equation}
Further properties of the joint modulus of variation are presented in Section~\ref{s:joint}.

\smallbreak
For a sequence of functions $\{f_j\}\subset M^T$ and $f\in M^T$, we write:
(a)~$f_j\to f$ on $T$ to denote the \emph{pointwise\/} (or \emph{everywhere\/})
\emph{convergence\/} of $\{f_j\}$ to $f$ (that is, $\lim_{j\to\infty}d(f_j(t),f(t))=0$
for all $t\in T$); (b) $f_j\rra f$ on $T$ to denote the \emph{uni\-form convergence\/}
of $\{f_j\}$ to $f$ meaning, as usual, that $\lim_{j\to\infty}d_\infty(f_j,f)=0$.
The uniform convergence implies the pointwise convergence, but not vice versa.
Recall that a sequence $\{f_j\}\subset M^T$ is said to be \emph{pointwise precompact\/}
on $T$ if the closure in $M$ of the set $\{f_j(t):j\in\NB\}$ is compact for all $t\in T$.

\smallbreak
Making use of E.~Landau's notation, given a sequence,
$\{\mu_n\}_{n=1}^\infty\subset\RB$, we write $\mu_n=o(n)$ to denote the
condition $\lim_{n\to\infty}\mu_n/n=0$.

Our main result, a \emph{pointwise selection principle\/} for metric space valued functions
in terms of the joint modulus of variation, is as follows.

\begin{thm} \label{t:first}
Let $\varnothing\ne T\subset\RB$ and $(M,d)$ be a metric space. Suppose
$\{f_j\},\{g_j\}\subset M^T$ are two sequences of functions such that\par\vspace{2pt}
{\rm(a)} $\{f_j\}$ is pointwise precompact on $T$,\par\vspace{2pt}
{\rm(b)} $\{g_j\}$ is pointwise convergent on $T$ to a function $g\in M^T$,
\par\vspace{2pt}\noindent
and
  \begin{equation} \label{e:limsup}
\mu_n\equiv\limsup_{j\to\infty}\nu_n(f_j,g_j)=o(n).
  \end{equation}
Then, there is a subsequence of $\{f_j\}$, which converges pointwise on $T$ to a
function $f\in M^T$ such that $\nu_n(f,g)\le\mu_n$ for all $n\in\NB$.
\end{thm}

This theorem will be proved in Section~\ref{s:proof}. Now, a few remarks are in order.
Given $f\in M^T$ and a constant function $c:T\to M$, the quantity
  \begin{equation} \label{e:nuf}
\nu_n(f)\equiv\nu_n(f,c)=\sup\biggl\{\sum_{i=1}^nd(f(s_i),f(t_i)):
\{I_i\}_1^n\prec T\biggr\}
  \end{equation}
(with $I_i=\{s_i,t_i\}$) is independent of $c$, and the sequence
\mbox{$\{\nu_n(f)\}_{n=1}^\infty\subset[0,\infty]$} is known as the
\emph{modulus of variation\/} of $f$ in the sense of Chanturiya (\cite{Chan74,Chan75,jmaa05,Ischia,MatTr,GNW}). It characterizes regulated
(or proper) functions on $T=[a,b]$ as follows. We say that $f:[a,b]\to M$ is
\emph{regulated\/} and write $f\in\mbox{\rm Reg}([a,b];M)$ if
$d(f(s),f(t))\to0$ as $s,t\to\tau-0$ for every $a<\tau\le b$, and
$d(f(s),f(t))\to0$ as $s,t\to\tau'+0$ for every $a\le\tau'<b$
(and so, by Cauchy's criterion, one-sided limits $f(\tau-0),f(\tau'+0)\in M$ exist
provided $M$ is complete). We have
  \begin{equation} \label{e:chan}
\mbox{\rm Reg}([a,b];M)=\{f\in M^{[a,b]}:\nu_n(f)=o(n)\}
  \end{equation}
(more general characterizations for dense subsets $T$ of $[a,b]$ can be found
in \cite{Ischia,MatTr}). A certain relationship between characterizations of regulated
functions and pointwise selection principles is exhibited in \cite{waterman80}.

\section{The joint modulus of variation} \label{s:joint}

We begin by studying the joint increment \eq{e:join}, whose properties are gathered
in the following lemma.

\begin{lem} \label{l:pjin}
Given $f,g,h\in M^T$ and $s,t\in T$, we have\/{\rm:}\par\vspace{2pt}
{\rm(a)} $|(f,f)(\{s,t\})|=0;$\par\vspace{2pt}
{\rm(b)} $|(f,g)(\{s,t\})|=|(g,f)(\{s,t\})|;$\par\vspace{2pt}
{\rm(c)} $|(f,g)(\{s,t\})|\le|(f,h)(\{s,t\})|+|(h,g)(\{s,t\})|;$\par\vspace{2pt}
{\rm(d)} $d(f(s),f(t))\le d(g(s),g(t))+|(f,g)(\{s,t\})|;$\par\vspace{2pt}
{\rm(e)} $d(f(t),g(t))\le d(f(s),g(s))+|(f,g)(\{s,t\})|$.
\end{lem}

\begin{pf}
Properties (a), (b), and (c), showing that $(f,g)\mapsto|(f,g)(\{s,t\})|$ is a
\emph{pseudometric\/} on $M^T$, are straightforward. To establish (d) and (e),
take into account equality $d(x,y)=\max_{z\in M}|d(x,z)-d(y,z)|$.
\qed\end{pf}

\begin{rem} \label{r:1}
(a) If $|(f,g)(\{s,t\})|=0$, then (d), (e), and (b) imply equalities
$d(f(s),f(t))=d(g(s),g(t))$ and $d(f(t),g(t))=d(f(s),g(s))$.
In addition to Lemma~\ref{l:pjin}, the function
$(s,t)\mapsto|(f,g)(\{s,t\})|$ is a pseudometric on $T$.

\smallbreak
(b) If $F(z)$ denotes the absolute value under the supremum sign in \eq{e:join}, then
$F:M\to[0,\infty)$ and $|F(z)-F(z_0)|\le4d(z,z_0)$ for all $z,z_0\in M$.

\smallbreak
(c) By Lemma~\ref{l:pjin}(d), $|f(T)|\le|g(T)|+|(f,g)(T)|=|g(T)|+\nu_1(f,g)$. So,
  $$\bigl||f(T)|-|g(T)|\bigr|\le|(f,g)(T)|\le|f(T)|+|g(T)|,\qquad f,g\in\mbox{\rm B}(T;M).$$
Moreover, it follows from Lemma~\ref{l:pjin}(e) that
  $$d_\infty(f,g)\le d(f(s),g(s))+|(f,g)(T)|\le3d_\infty(f,g)\quad\,
     \mbox{for \,all}\quad\,s\in T.$$

(d) Suppose the triple $(M,d,+)$ is a \emph{metric semigroup\/} (\cite[Section~4]{Sovae}),
i.e., $(M,d)$ is a metric space, $(M,+)$ is an Abelian semigroup with the operation of
addition~$+$, and $d(x,y)=d(x+z,y+z)$ for all $x,y,z\in M$. Then, the joint increment
\eq{e:join} may be alternatively replaced by
  \begin{equation} \label{e:mes}
|(f,g)(\{s,t\})|=d(f(s)+g(t),f(t)+g(s)).
  \end{equation}
The joint modulus of variation \eq{e:jmv} involving \eq{e:mes} was employed
in \cite{jmaa08}. Furthermore, if $(M,\|\cdot\|)$ is a normed linear space
(over $\RB$ or $\mathbb{C}$), we may set
  \begin{equation} \label{e:nrm}
|(f,g)(\{s,t\})|=\|f(s)+g(t)-f(t)-g(s)\|=\|(f-g)(s)-(f-g)(t)\|.
  \end{equation}
Quantities \eq{e:mes} and \eq{e:nrm} have the same properties as \eq{e:join}:
see \eq{e:fg1}, \eq{e:fg2}, Lemma~\ref{l:pjin} and Remark~\ref{r:1}(a).
In the sequel, we make use of more general quantity \eq{e:join}.
\end{rem}

If $f,g\in M^T$, $n\in\NB$ and $\varnothing\ne E\subset T$, we set
$\nu_n(f,g;E)=\nu_n(f|_E,g|_E)$, where $f|_E\in M^E$ is the restriction of $f$ to~$E$
(i.e., $f|_E(t)=f(t)$ for all $t\in E$). Accordingly, $\nu_n(f,g)=\nu_n(f,g;T)$.

\smallbreak
The following properties of the joint modulus of variation are immediate.
The sequence $\{\nu_n(f,g)\}_{n=1}^\infty$ is nondecreasing,
$\nu_{n+m}(f,g)\le\nu_n(f,g)+\nu_m(f,g)$ for all $n,m\in\NB$, and
$\nu_n(f,g;E)\le\nu_n(f,g;T)$ provided $n\in\NB$ and $E\subset T$.
It follows from \eq{e:jmv} and Lemma~\ref{l:pjin}(a)--(c) that, for every $n\in\NB$,
the function $(f,g)\mapsto\nu_n(f,g)$ is a \emph{pseudometric\/} on $M^T$
(possibly assuming infinite values) and, in particular (cf.\ \eq{e:jmv} and \eq{e:nuf}),
  \begin{equation} \label{e:pp}
\nu_n(f,g)\le\nu_n(f)+\nu_n(g)\quad\mbox{ and }\quad
\nu_n(f)\le\nu_n(g)+\nu_n(f,g)
  \end{equation}
for all $n\in\NB$ and $f,g\in M^T$. Furthemore, if $f,g\in\mbox{\rm B}(T;M)$, then, by 
\eq{e:nuon}, the sequence $\{\nu_n(f,g)/n\}_{n=1}^\infty$ is bounded in $[0,\infty)$.

\pagebreak
Essential properties of the joint modulus of variation are presented in

\begin{lem} \label{l:ess}
Given $n\in\NB$, $f,g\in M^T$, and $\varnothing\ne E\subset T$, we have\/{\rm:}
\par\vspace{2pt}
{\rm(a)} $|(f,g)(\{s,t\})|+\nu_n(f,g;E_s^-)\le\nu_{n+1}(f,g;E_t^-)$ for all\/ $s,t\in E$
  with\par\vspace{2pt}\quad\,\,
  $s\le t$, where $E_\tau^-=(-\infty,\tau]\cap E$ for $\tau\in E;$
\par\vspace{2pt}
{\rm(b)} $\displaystyle\nu_{n+1}(f,g;E)\le\nu_n(f,g;E)+\frac{\nu_{n+1}(f,g;E)}{n+1}\,;$
\par\vspace{2pt}
{\rm(c)} if\/ $\{f_j\},\{g_j\}\subset M^T$ are such that $f_j\to f$
  and $g_j\to g$ on $E$, then\par\vspace{2pt}\quad\,\,
  $\nu_n(f,g;E)\le\liminf_{j\to\infty}\nu_n(f_j,g_j;E);$
\par\vspace{2pt}
{\rm(d)} if\/ $\{f_j\},\{g_j\}\subset M^T$ are such that $f_j\rra f$
  and $g_j\rra g$ on $E$, then\par\vspace{2pt}\quad\,\,
  $\nu_n(f,g;E)=\lim_{j\to\infty}\nu_n(f_j,g_j;E).$
\end{lem}

\begin{pf}
(a) We may assume that $s<t$. Let $\{I_i\}_1^n\prec E_s^-$. Setting $I_0=\{s,t\}$,
we find $\{I_i\}_0^n\prec E_t^-$, and so,
  $$|(f,g)(I_0)|+\sum_{i=1}^n|(f,g)(I_i)|\le\nu_{n+1}(f,g;E_t^-).$$
The inequality in (a) follows by taking the supremum over all $\{I_i\}_1^n\prec E_s^-$.

\smallbreak
(b) We may assume that $\nu_{n+1}(f,g;E)$ is finite, and apply the idea from
\cite[Lemma]{Chan75}. Given $\vep>0$, there is $\{I_i\}_1^{n+1}\prec E$
(depending on~$\vep$) such that
  $$\sum_{i=1}^{n+1}|(f,g)(I_i)|\le\nu_{n+1}(f,g;E)
     \le\sum_{i=1}^{n+1}|(f,g)(I_i)|+\vep.$$
Setting $a_0=\min_{1\le i\le n+1}|(f,g)(I_i)|$, the left-hand side inequality implies
\mbox{$(n+1)a_0\le\nu_{n+1}(f,g;E)$}. The right-hand side inequality gives
  $$\nu_{n+1}(f,g;E)\le\nu_n(f,g;E)+a_0+\vep,$$
from which our inequality follows due to the arbitrariness of $\vep>0$.

\smallbreak
(c) First, we note that, given $j\in\NB$ and $s,t\in T$, we have
  \begin{align}
\bigl||(f_j,g_j)(\{s,t\})|-|(f,g)(\{s,t\})|\bigr|&\le d(f_j(s),f(s))+d(f_j(t),f(t))
  \nonumber\\[2pt]
&\quad\,\,+d(g_j(s),g(s))+d(g_j(t),g(t)). \label{e:ten}
  \end{align}
In fact, Lemma~\ref{l:pjin}(c) and inequality \eq{e:fg2} imply
  \begin{align}
|(f_j,g_j)(\{s,t\})|&\le|(f_j,f)(\{s,t\})|+|(f,g)(\{s,t\})|+|(g,g_j)(\{s,t\})|
  \nonumber\\[2pt]
&\le d(f_j(s),f(s))+d(f_j(t),f(t))+|(f,g)(\{s,t\})|\nonumber\\[2pt]
&\quad\,\,+d(g(s),g_j(s))+d(g(t),g_j(t)). \label{e:sstt}
  \end{align}
Exchanging $f_j$ and $f$ as well as $g_j$ and $g$, we obtain \eq{e:ten}.

\smallbreak
From the pointwise convergence of $\{f_j\}$ and $\{g_j\}$ and \eq{e:ten}, we find
  $$\lim_{j\to\infty}|(f_j,g_j)(\{s,t\})|=|(f,g)(\{s,t\})|\quad\,
     \mbox{for \,all}\quad\,s,t\in E.$$
By definition \eq{e:jmv}, given $\{I_i\}_1^n\prec E$, we have
  $$\sum_{i=1}^n|(f_j,g_j)(I_i)|\le\nu_n(f_j,g_j;E)\quad\,
     \mbox{for \,all}\quad\,j\in\NB.$$
Passing to the limit inferior as $j\to\infty$, we get
  \begin{equation} \label{e:linf}
\sum_{i=1}^n|(f,g)(I_i)|\le\liminf_{j\to\infty}\nu_n(f_j,g_j;E).
  \end{equation}
Since $\{I_i\}_1^n\prec E$ is arbitrary, it remains to take into account \eq{e:jmv}.

\smallbreak
(d) It follows from \eq{e:sstt} that, for any $s,t\in E$ and $j\in\NB$,
  $$|(f_j,g_j)(\{s,t\})|\!\le\!2\sup_{\tau\in E}d(f_j(\tau),f(\tau))\!+\!|(f,g)(\{s,t\})
     \!+\!2\sup_{\tau\in E}d(g_j(\tau),g(\tau)),$$
and so, definition \eq{e:jmv} implies
  \begin{equation} \label{e:Et}
\nu_n(f_j,g_j;E)\!\le\!2n\sup_{\tau\in E}d(f_j(\tau),f(\tau))\!+\!\nu_n(f,g;E)
\!+\!2n\sup_{\tau\in E}d(g_j(\tau),g(\tau))
  \end{equation}
for all $j\in\NB$. Passing to the limit superior as $j\to\infty$, we get
  $$\limsup_{j\to\infty}\nu_n(f_j,g_j;E)\le\nu_n(f,g;E).$$
Now, the equality in (d) is a consequence of Lemma~\ref{l:ess}(c).
\qed\end{pf}

\begin{rem} \label{r:2}
If the value $\nu_1(f,g;E)=|(f,g)(E)|$ (see \eq{e:nuon}) is finite for an $E\subset T$
(e.g., when $f,g\in\mbox{\rm B}(E;M)$), inequality in Lemma~\ref{l:ess}(b)
is equivalent to
  $$\frac{\nu_{n+1}(f,g;E)}{n+1}\le\frac{\nu_n(f,g;E)}{n}\,.$$
Thus, the limit $\lim_{n\to\infty}\nu_n(f,g;E)/n$ always exists in $[0,\infty)$.
\end{rem}

\section{Conditionally regulated functions} \label{s:reg}

Since $\nu_n=\nu_n(\cdot,\cdot)$ is a (extended) pseudometric on $M^T$, we may
introduce an equivalence relation $\sim$ on $M^T$ as follows: given $f,g\in M^T$, we set
  $$\mbox{$f\sim g$ \,\,\,if \,and \,only \,if \,\,\,$\nu_n(f,g)=o(n)$.}$$
The equivalence class $\mbox{\rm R}(g)=\{f\in M^T:f\sim g\}$ of a function $g\in M^T$
is called the \emph{regularity class\/} of $g$, and any representative
$f\in\mbox{\rm R}(g)$ is called a \emph{conditionally regulated\/} or, more precisely,
\emph{$g$-regulated\/} function. This terminology is justified by \eq{e:chan}: in the
framework of the product set $M^{[a,b]}$, we have $\mbox{Reg}([a,b];M)=\mbox{R}(c)$
for any constant function~$c:[a,b]\to M$.

\smallbreak
Note that, in Theorem~\ref{t:first}, condition `$\nu_n(f,g)\le\mu_n$ for all $n\in\NB$'
means that $f\in\mbox{\rm R}(g)$, and so, the class $\mbox{\rm R}(g)$ is worth
studying in more detail.

\begin{thm} \label{t:greg}
Given a function $g\in M^T$, we have\/{\rm:}\par\vspace{2pt}
{\rm(a)} $g\in\mbox{\rm B}(T;M)$ if and only if\/
  $\mbox{\rm R}(g)\subset\mbox{\rm B}(T;M);$\par\vspace{2pt}
{\rm(b)} $\mbox{\rm R}(g)$ is closed with respect to the uniform convergence,
  but not closed\par\quad\,\,
 with respect to the pointwise convergence in general\/{\rm;}
\par\vspace{2pt}
{\rm(c)} if $(M,d)$ is a complete metric space, then the pair
  $(\mbox{\rm R}(g),d_\infty)$ is also\par\quad\,\,
   a complete metric space.
\end{thm}

\begin{pf}
(a) The sufficiency is clear, because $g\in\mbox{\rm R}(g)$. Now, suppose that
$g\in\mbox{\rm B}(T;M)$, so that, by \eq{e:nuf}, $\nu_1(g)=|(g,c)(T)|=|g(T)|<\infty$.
Given $f\in\mbox{\rm R}(g)$, $\nu_n(f,g)=o(n)$, and so, $\nu_{n_0}(f,g)\le n_0$
for some $n_0\in\NB$. It follows from \eq{e:pp} and \eq{e:nuon} that
  $$|f(T)|=\nu_1(f)\le\nu_1(g)+\nu_1(f,g)\le|g(T)|+\nu_{n_0}(f,g)\le|g(T)|+n_0<\infty,$$
which implies $f\in\mbox{\rm B}(T;M)$.

\smallbreak
(b) It is to be shown that if $\{f_j\}\subset\mbox{\rm R}(g)$ and $f_j\rra f$ on $T$
with $f\in M^T$, then $f\in\mbox{\rm R}(g)$. We will prove a little bit more: namely,
if $\{f_j\},\{g_j\}\subset M^T$ are such that $f_j\in\mbox{\rm R}(g_j)$ for all $j\in\NB$,
$f_j\rra f$ and $g_j\rra g$ on $T$ with $f,g\in M^T$, then $f\in\mbox{\rm R}(g)$
(the previous assertion follows if $g_j=g$ for all $j\in\NB$).
In fact, exchanging $f_j$ and $f$, and $g_j$ and $g$ in \eq{e:Et}, we get
  $$\frac{\nu_n(f,g)}n\le2d_\infty(f,f_j)+\frac{\nu_n(f_j,g_j)}n+2d_\infty(g,g_j),
     \quad\,\,n,j\in\NB.$$
By the uniform convergence of $\{f_j\}$ and $\{g_j\}$, given $\vep>0$, there is
a number $j_0=j_0(\vep)\in\NB$ such that $d_\infty(f,f_{j_0})\le\vep$ and
$d_\infty(g,g_{j_0})\le\vep$. Since $f_{j_0}$ is in $\mbox{\rm R}(g_{j_0})$, we have
$\nu_n(f_{j_0},g_{j_0})=o(n)$, and so, there exists $n_0=n_0(\vep)\in\NB$ such that
$\nu_n(f_{j_0},g_{j_0})/n\le\vep$ for all $n\ge n_0$. The estimate above with $j=j_0$
implies $\nu_n(f,g)/n\le5\vep$, $n\ge n_0$, which means that
$\nu_n(f,g)=o(n)$ and $f\in\mbox{\rm R}(g)$.

\smallbreak
As for the pointwise convergence, consider a sequence of real step functions converging
pointwise to the Dirichlet function ($=$\,the characteristic function of the rationals~%
$\mathbb{Q}$) on $T=[0,1]$ (see \cite[Examples~4, 5]{jmaa05} and
Example~\ref{ex:Di}(a) in Section~\ref{s:proof}).

\smallbreak
(c) First, we show that $d_\infty(f,f')<\infty$ for all $f,f'\in\mbox{\rm R}(g)$. In fact,
since $f\sim f'$, we have $\nu_n(f,f')=o(n)$, and so, $\nu_{n_0}(f,f')\le n_0$ for
some $n_0\in\NB$. Given $s\in T$, it follows from Remark~\ref{r:1}(c) and
\eq{e:nuon} that
  $$d_\infty(f,f')\le d(f(s),f'(s))+\nu_1(f,f')\le d(f(s),f'(s))+\nu_{n_0}(f,f')<\infty.$$

\smallbreak
The metric axioms for $d_\infty$ on $\mbox{\rm R}(g)$ are verified in a standard way.

\smallbreak
In order to prove that $\mbox{\rm R}(g)$ is complete, suppose
$\{f_j\}\subset\mbox{\rm R}(g)$ is a Cauchy sequence, i.e., $d_\infty(f_j,f_k)\to0$
as $j,k\to\infty$. Since $d(f_j(t),f_k(t))\le d_\infty(f_j,f_k)$ for all $t\in T$ and
$(M,d)$ is complete, there exists $f\in M^T$ such that $f_j\to f$ on $T$.
Noting that $f_j\to f_j$ and $f_k\to f$ on $T$ as $k\to\infty$ (and arguing as in
\eq{e:linf}), we get
  $$d_\infty(f_j,f)\le\liminf_{k\to\infty}d_\infty(f_j,f_k)
     =\lim_{k\to\infty}d_\infty(f_j,f_k)<\infty\quad\mbox{for all}\quad j\in\NB.$$
Since the sequence $\{f_j\}$ is $d_\infty$-Cauchy, we find
  $$\limsup_{j\to\infty}d_\infty(f_j,f)\le
     \lim_{j\to\infty}\lim_{k\to\infty}d_\infty(f_j,f_k)=0.$$
Thus, $\lim_{j\to\infty}d_\infty(f_j,f)=0$, and so, $f_j\rra f$ on $T$. Applying item (b)
of this Theorem, we conclude that $f\in\mbox{\rm R}(g)$.
\qed\end{pf}

A traditionally important class of regulated functions is the space of functions of bounded
Jordan variation, $\mbox{\rm BV}(T;M)$, which is introduced by means of the joint
modulus of variation as follows.

\smallbreak
Since the sequence $\{\nu_n(f,g)\}_{n=1}^\infty$ is nondecreasing for all
$f,g\in M^T$, the quantity (finite or not)
  $V(f,g)=\lim_{n\to\infty}\nu_n(f,g)=\sup_{n\in\NB}\nu_n(f,g)$
is called the \emph{joint variation\/} of functions $f$ and $g$ on $T$. The value
$V(f)\equiv V(f,c)$ is independent of a constant function $c:T\to M$ and is the
usual Jordan \emph{variation\/} of $f$ on $T$:
  $$V(f)=\sup\biggl\{\sum_{i=1}^nd(f(t_i),f(t_{i-1})):\mbox{$n\in\NB$ and
      $\{t_i\}_0^n\prec T$}\biggr\},$$
the supremum being taken over all partitions $\{t_i\}_0^n$ of $T$ (cf.\
Section~\ref{s:main}). The set $\mbox{\rm BV}(T;M)=\{f\in M^T:V(f)<\infty\}$
is contained in $\mbox{\rm B}(T;M)\cap\mbox{\rm R}(c)$
(in fact, $|f(T)|=\nu_1(f)\le V(f)$ and $\nu_n(f,c)/n\le V(f)/n$ for all
$f\in\mbox{\rm BV}(T;M)$).

\smallbreak
The following notion of \emph{$\vep$-variation\/} $V_\vep(f)$, due to Fra{\v n}kov{\'a}
\cite[Section~3]{Fr}, provides an alternative characterization (cf.\ \eq{e:chan}) of
regulated functions: given $f\in M^T$ and $\vep>0$, set
  \begin{equation} \label{e:evar}
V_\vep(f)\!=\!\inf\,\bigl\{V(g):\mbox{$g\in\mbox{\rm BV}(T;M)$ and
$d_\infty(f,g)\le\vep$}\bigr\}\quad(\inf\varnothing\!=\!\infty).
  \end{equation}
It was shown in \cite[Proposition~3.4]{Fr} (for $T=[a,b]$) that
  \begin{equation} \label{e:fra}
\mbox{\rm Reg}([a,b];M)=\{f\in M^{[a,b]}:\mbox{$V_\vep(f)<\infty$ for all $\vep>0$}\}
  \end{equation}
(although it was assumed in \cite{Fr} that $M=\RB^N$, the proof of the last assertion
carries over to any metric space $M$, cf.\ \cite[Lemma~3]{jmaa05}).

\smallbreak
The notion of $\vep$-variation will be needed in Section~\ref{s:ext}.

\begin{exa} \label{ex:Dir-type}
Given $x,y\!\in\! M$ with $x\!\ne\! y$, let $f\!=\!\DC_{x,y}:T\!=\![0,1]\to M$
be the Dirichlet-type function of the form:
  \begin{equation} \label{e:Dxy}
\DC_{x,y}(t)=\left\{
   \begin{array}{ccl}
x & \mbox{if} & \mbox{$t\in[0,1]$ is rational,}\\[2pt]
y & \mbox{if} & \mbox{$t\in[0,1]$ is irrational.}
   \end{array}\right.
  \end{equation}
Clearly, $f\notin\mbox{\rm Reg}([0,1];M)$. Moreover (cf.\ \eq{e:fra}), we have
  \begin{equation} \label{e:Vef}
\mbox{$V_\vep(f)=\infty$ \,if \,$0<\vep<d(x,y)/2$,\quad and \,
$V_\vep(f)=0$ \,if \,$\vep\ge d(x,y)$.}
  \end{equation}
To see this, first note that, given $g\in M^{[0,1]}$, inequality $d_\infty(f,g)\le\vep$
is equivalent to the following two conditions:
  \begin{equation} \label{e:Quf}
\mbox{$d(x,g(s))\le\vep$ \,$\forall\,s\in[0,1]\cap\mathbb{Q}$,\quad and \,
$d(y,g(t))\le\vep$ \,$\forall\,t\in[0,1]\setminus\mathbb{Q}$.}
  \end{equation}

Suppose $0<\vep<d(x,y)/2$. To show that $d_\infty(f,g)\le\vep$ implies $V(g)=\infty$,
we let $n\in\NB$, and $\{t_i\}_0^{2n}\prec[0,1]$ be a partition of $[0,1]$ such that
points  $\{t_{2i}\}_{i=0}^n$ are rational and points $\{t_{2i-1}\}_{i=1}^n$ are irrational.
By the triangle inequality for $d$ and \eq{e:Quf}, we get
  \begin{align*}
V(g)&\ge\sum_{i=1}^{2n}d(g(t_i),g(t_{i-1}))\ge\sum_{i=1}^nd(g(t_{2i}),g(t_{2i-1}))\\
&\ge\sum_{i=1}^n\bigl(d(x,y)-d(x,g(t_{2i}))-d(g(t_{2i-1}),y)\bigr)
  \ge n\bigl(d(x,y)-2\vep\bigr).
  \end{align*}

If $\vep\ge d(x,y)$, we set $g(t)=x$ (or $g(t)=y$) for all $t\in[0,1]$, so that
\eq{e:Quf} is satisfied and $V(g)=0$. Thus, $V_\vep(f)=0$.

\smallbreak
The second assertion in \eq{e:Vef} can be refined, provided
  \begin{equation} \label{e:vex}
d(x,y)/2=\max\{d(x,z_0),d(y,z_0)\}\quad\mbox{for \,some}\quad z_0\in M.
  \end{equation}
In fact, we may set $g(t)=z_0$ for all $t\in[0,1]$, so that \eq{e:Quf} holds whenever
$d(x,y)/2\le\vep$, and
$V(g)=0$. This implies $V_\vep(f)=0$ for all $\vep\ge d(x,y)/2$.

\smallbreak
A few remarks concerning condition \eq{e:vex} are in order. Since
  $$d(x,y)\le d(x,z)+d(z,y)\le2\max\{d(x,z),d(y,z)\}\quad\mbox{for \,all}\quad z\in M,$$
condition \eq{e:vex} refers to a certain form of `convexity' of $M$ (which is not
restrictive for our purposes). For instance, if $(M,\|\cdot\|)$ is a normed linear space
with $d(x,y)=\|x-y\|$, we may set $z_0=(x+y)/2$. More generally, by Menger's
Theorem (\cite{Menger}, \cite[Example 2.7]{Go-Ki}), if a metric space $(M,d)$ is
complete and \emph{metrically convex\/} (i.e., given $x,y\in M$ with $x\ne y$, there is
$z\in M$ such that $x\ne z\ne y$ and $d(x,y)=d(x,z)+d(z,y)$), then, for any
$x,y\in M$, there is an isometry $\varphi:[0,d(x,y)]\to M$ such that $\varphi(0)=x$
and $\varphi(d(x,y))=y$. In this case, we set $z_0=\varphi(d(x,y)/2)$.
More examples of metrically convex metric spaces can be found in \cite{ChR,Duda}.

\smallbreak
Finally, if $M=\{x,y\}$, then condition \eq{e:vex} is not satisfied, and we have
$V_\vep(f)=\infty$ for all $0<\vep<d(x,y)$, which is a consequence of \eq{e:Quf}.
\end{exa}

\section{Proof of the main result} \label{s:proof}

\begin{pfirst}
With no loss of generality we may assume that $T$ is uncountable; otherwise, by virtue
of assumption (a) and the standard Cantor diagonal procedure, we extract a pointwise
convergent subsequence of $\{f_j\}$ and apply Lemma~\ref{l:ess}(c). Note that $\mu_n$
is finite for all $n\in\NB$: in fact, $\mu_n\le n$ whenever $n\ge n_0$ for some $n_0\in\NB$
and, since $n\mapsto\nu_n(f_j,g_j)$ is nondecreasing for all $j\in\NB$, we have
$\mu_n\le\mu_{n_0}\le n_0$ for all $1\le n<n_0$.

\smallbreak
We divide the rest of the proof into four steps for clarity.

\smallbreak
\emph{Step~1.} Let us show that there is a subsequence of $\{j\}_{j=1}^\infty$,
again denoted by $\{j\}$, and a nondecreasing sequence
$\{\alpha_n\}_{n=1}^\infty\subset[0,\infty)$ such that
  \begin{equation} \label{e:alf}
\lim_{j\to\infty}\nu_n(f_j,g_j)=\alpha_n\le\mu_n\quad\,\mbox{for \,all}\quad\,n\in\NB.
  \end{equation}

\smallbreak
We set $\alpha_1=\mu_1$. The definition \eq{e:limsup} of $\mu_1$ implies that
there is an increasing sequence $\{J_1(j)\}_{j=1}^\infty\subset\NB$ (i.e., a
subsequence of $\{j\}_{j=1}^\infty$) such that
$\nu_1(f_{J_1(j)},g_{J_1(j)})\to\alpha_1$ as $j\to\infty$. Setting
$\alpha_2=\limsup_{j\to\infty}\nu_2(f_{J_1(j)},g_{J_1(j)})$, we find $\alpha_2\le\mu_2$,
and there is a subsequence $\{J_2(j)\}_{j=1}^\infty$ of $\{J_1(j)\}_{j=1}^\infty$
such that $\nu_2(f_{J_2(j)},g_{J_2(j)})\to\alpha_2$ as $j\to\infty$. Inductively,
if $n\ge3$ and a subsequence $\{J_{n-1}(j)\}_{j=1}^\infty$ of $\{j\}_{j=1}^\infty$
is already chosen, we define $\alpha_n$ as the limit superior of
$\nu_n(f_{J_{n-1}(j)},g_{J_{n-1}(j)})$ as $j\to\infty$,
so that $\alpha_n\le\mu_n$. Now, we pick a subsequence $\{J_n(j)\}_{j=1}^\infty$ of
$\{J_{n-1}(j)\}_{j=1}^\infty$ such that $\nu_n(f_{J_n(j)},g_{J_n(j)})\to\alpha_n$
as $j\to\infty$. Noting that the sequence $\{J_j(j)\}_{j=n}^\infty$ is a subsequence
of $\{J_n(j)\}_{j=1}^\infty$ (for all $n\in\NB$) and denoting the diagonal
sequences $\{f_{J_j(j)}\}_{j=1}^\infty$ and $\{g_{J_j(j)}\}_{j=1}^\infty$ again by
$\{f_j\}$ and $\{g_j\}$, respectively, we obtain \eq{e:alf}.

\smallbreak
In the sequel, the set of all nondecreasing bounded functions mapping $T$ into
$\RB^+=[0,\infty)$ is denoted by $\mbox{\rm Mon}(T;\RB^+)$. 

\smallbreak
\emph{Step~2.} In this step, we prove that there are subsequences of $\{f_j\}$ and
$\{g_j\}$ from \eq{e:alf}, again denoted by $\{f_j\}$ and $\{g_j\}$, respectively,
and a sequence of functions $\{\beta_n\}_{n=1}^\infty\subset\mbox{Mon}(T;\RB^+)$
such that
  \begin{equation} \label{e:bet}
\lim_{j\to\infty}\nu_n(f_j,g_j;T_t^-)=\beta_n(t)\quad\mbox{for \,all \,$n\in\NB$
\,and \,$t\in T$,}
  \end{equation}
where $T_t^-=\{s\in T:s\le t\}$ for $t\in T$.

\smallbreak
Note that, for each $n\in\NB$, the function $t\mapsto\nu_n(f_j,g_j;T_t^-)$ is
nondecreasing on $T$, and $\nu_n(f_j,g_j;T_t^-)\le\nu_n(f_j,g_j)$ for all $t\in T$
and $n\in\NB$. By virtue of \eq{e:alf}, there is a sequence
$\{C_n\}_{n=1}^\infty\subset\RB^+$ such that $\nu_n(f_j,g_j)\le C_n$ for all
$n,j\in\NB$. In what follows, we apply the diagonal procedure once again.

\smallbreak
The sequence $\{t\mapsto\nu_1(f_j,g_j;T_t^-)\}_{j=1}^\infty%
\subset\mbox{\rm Mon}(T;\RB^+)$ is uniformly bounded by constant $C_1$, and so,
by Helly's Selection Principle, there are an increasing sequence
$\{K_1(j)\}_{j=1}^\infty\subset\NB$ (i.e., a subsequence of $\{j\}_{j=1}^\infty$)
and a function $\beta_1\in\mbox{\rm Mon}(T;\RB^+)$ such that
$\nu_1(f_{K_1(j)},g_{K_1(j)};T_t^-)\to\beta_1(t)$ as $j\to\infty$ for all $t\in T$.
The sequence $\{t\mapsto\nu_2(f_{K_1(j)},g_{K_1(j)};T_t^-)\}_{j=1}^\infty%
\subset\mbox{\rm Mon}(T;\RB^+)$ is uniformly bounded on $T$ by constant $C_2$,
and so, again by Helly's Theorem, there are a subsequence
$\{K_2(j)\}_{j=1}^\infty$ of $\{K_1(j)\}_{j=1}^\infty$ and a function
$\beta_2\in\mbox{\rm Mon}(T;\RB^+)$ such that
$\nu_2(f_{K_2(j)},g_{K_2(j)};T_t^-)\to\beta_2(t)$ as $j\to\infty$ for all $t\in T$.
Inductively, if $n\ge3$ and a subsequence $\{K_{n-1}(j)\}_{j=1}^\infty$ of
$\{j\}_{j=1}^\infty$ and a function $\beta_{n-1}\in\mbox{\rm Mon}(T;\RB^+)$
are already chosen, we apply the Helly Theorem to the sequence of functions
$\{t\mapsto\nu_n(f_{K_{n-1}(j)},g_{K_{n-1}(j)};T_t^-)\}_{j=1}^\infty%
\subset\mbox{\rm Mon}(T;\RB^+)$, which is uniformly bounded on $T$
by constant $C_n$: there are a subsequence $\{K_n(j)\}_{j=1}^\infty$ of
$\{K_{n-1}(j)\}_{j=1}^\infty$ and a function $\beta_n\in\mbox{\rm Mon}(T;\RB^+)$
such that $\nu_n(f_{K_n(j)},g_{K_n(j)};T_t^-)\to\beta_n(t)$ as $j\to\infty$ for
all $t\in T$. Since $\{K_j(j)\}_{j=n}^\infty$ is a subsequence of
$\{K_n(j)\}_{j=1}^\infty$ (for all $n\in\NB$), the diagonal sequences
$\{f_{K_j(j)}\}_{j=1}^\infty$ and $\{g_{K_j(j)}\}_{j=1}^\infty$, again denoted by
$\{f_j\}$ and $\{g_j\}$, respectively, satisfy condition~\eq{e:bet}.

\smallbreak
\emph{Step~3.} Let $Q$ be an at most countable dense subset of $T$.
Note that~$Q$ contains all isolated ($=$\,non-limit) points of $T$ (i.e., points $t\in T$
such that the intervals $(t-\delta,t)$ and $(t,t+\delta)$ lie in $\RB\setminus T$ for
some $\delta>0$). The set $Q_n\subset T$ of discontinuity points of nondecreasing
function $\beta_n$ is at most countable. Setting $S=Q\cup\bigcup_{n=1}^\infty Q_n$,
we find that $S$ is an at most countable dense subset of $T$ and
  \begin{equation} \label{e:Set}
\mbox{$\beta_n$ is continuous at all points of $T\setminus S$ for all $n\in\NB$.}
  \end{equation}
Since the set $\{f_j(t):j\in\NB\}$ is precompact in $M$ for all $t\in T$, and
$S\subset T$ is at most countable, we may assume (applying the diagonal procedure
again and passing to a subsequence of $\{f_j\}$ if necessary) that, given $s\in S$,
there is a point $f(s)\in M$ such that $d(f_j(s),f(s))\to0$ as $j\to\infty$. In this way,
we obtain a function $f:S\to M$.

\smallbreak
\emph{Step~4.} Now, we show that, for every $t\in T\setminus S$, the sequence
$\{f_j(t)\}_{j=1}^\infty$ converges in $M$. For this, we prove that this sequence is
Cauchy in $M$, i.e., $d(f_j(t),f_k(t))\to0$ as $j,k\to\infty$. Let $\vep>0$ be arbitrarily
fixed. By assumption \eq{e:limsup}, $\mu_n/n\to0$ as $n\to\infty$, so we choose and
fix a number $n=n(\vep)\in\NB$ such that
  $$\frac{\mu_{n+1}}{n+1}\le\vep.$$
By property \eq{e:alf}, there is a number $j_1=j_1(\vep,n)\in\NB$ such that
  $$\nu_{n+1}(f_j,g_j)\le\alpha_{n+1}+\vep\le\mu_{n+1}+\vep\quad
     \mbox{for \,all}\quad j\ge j_1.$$
The definition of the set $S$ and \eq{e:Set} imply that the point $t$ is a limit point
for $T$ and, at the same time, a point of continuity of the function $\beta_n$.
By the density of $S$ in $T$, there is $s=s(\vep,n,t)\in S$ such that
  $$|\beta_n(t)-\beta_n(s)|\le\vep.$$
It follows from \eq{e:bet} that there is $j_2=j_2(\vep,n,t,s)\in\NB$ such that
  $$|\nu_n(f_j,g_j;T_t^-)-\beta_n(t)|\le\vep\quad\mbox{and}\quad
     |\nu_n(f_j,g_j;T_s^-)-\beta_n(s)|\le\vep\quad\forall\,j\ge j_2.$$
Assuming that $s<t$ (the arguments are similar if $t<s$) and applying
Lemma~\ref{l:ess}(a),\,(b), we get, for all $j\ge\max\{j_1,j_2\}$:
  \begin{align*}
|(f_j,g_j)(\{s,t\})|&\le\nu_{n+1}(f_j,g_j;T_t^-)-\nu_n(f_j,g_j;T_s^-)\\[2pt]
&\le\nu_{n+1}(f_j,g_j;T_t^-)-\nu_n(f_j,g_j;T_t^-)\\[2pt]
&\quad\,+|\nu_n(f_j,g_j;T_t^-)-\beta_n(t)|+|\beta_n(t)-\beta_n(s)|\\[2pt]
&\quad\,+|\beta_n(s)-\nu_n(f_j,g_j;T_s^-)|\\
&\le\frac{\nu_{n+1}(f_j,g_j;T_t^-)}{n+1}+\vep+\vep+\vep\\
&\le\frac{\nu_{n+1}(f_j,g_j)}{n+1}+3\vep\\
&\le\frac{\mu_{n+1}}{n+1}+\frac\vep{n+1}+3\vep\le5\vep.
  \end{align*}
Being convergent (see condition (b)), sequences $\{f_j(s)\}_{j=1}^\infty$,
$\{g_j(s)\}_{j=1}^\infty$ and $\{g_j(t)\}_{j=1}^\infty$ are Cauchy in $M$, and so,
there is a number $j_3=j_3(\vep,s,t)\in\NB$ such that, for all $j,k\ge j_3$, we have
  $$d(f_j(s),f_k(s))\le\vep,\quad d(g_j(s),g_k(s))\le\vep,\quad\mbox{and}\quad
     d(g_j(t),g_k(t))\le\vep.$$
By virtue of inequality \eq{e:fg2}, we get
  $$|(g_j,g_k)(\{s,t\})|\le d(g_j(s),g_k(s))+d(g_j(t),g_k(t))\le2\vep\quad
     \forall\,j,k\ge j_3.$$
Putting $j_4=\max\{j_1,j_2,j_3\}$ and applying Lemma~\ref{l:pjin}(e),\,(c),\,(b),
we find
  \begin{align*}
d(f_j(t),f_k(t))&\le d(f_j(s),f_k(s))+|(f_j,f_k)(\{s,t\})|\\[2pt]
&\le d(f_j(s),f_k(s))+|(f_j,g_j)(\{s,t\})|+|(g_j,g_k)(\{s,t\})|\\
&\quad\,+|(g_k,f_k)(\{s,t\})|\\[2pt]
&\le\vep+5\vep+2\vep+5\vep=13\vep\quad\mbox{for all}\quad j,k\ge j_4.
  \end{align*}
Since $j_4$ depends only on $\vep$ (and $t$), this proves that $\{f_j(t)\}_{j=1}^\infty$
is a Cauchy sequence in $M$, which together with assumption (a) establishes its
convergence in $M$ to an element denoted by $f(t)\in M$.

\smallbreak
Here and at the end of Step~3, we have shown that $f:T=S\cup(T\setminus S)\to M$
is a pointwise limit on $T$ of a subsequence $\{f_{j_k}\}_{k=1}^\infty$ of the original
sequence $\{f_j\}_{j=1}^\infty$. Since $g_{j_k}\to g$ pointwise on $T$ as
$k\to\infty$ as well, we conclude from Lemma~\ref{l:ess}(c) that
  $$\nu_n(f,g)\le\liminf_{k\to\infty}\nu_n(f_{j_k},g_{j_k})\le
     \limsup_{j\to\infty}\nu_n(f_j,g_j)=\mu_n\quad\forall\,n\in\NB,$$
and so, $\nu_n(f,g)\!=\!o(n)$, or $f\!\in\!\mbox{\rm R}(g)$.
This completes the proof of Theorem~\ref{t:first}.\,\,$\Box$
\end{pfirst}

\begin{rem}
(a) Condition (b) in Theorem~\ref{t:first} may be replaced by the following one:
$\{g_j(t)\}_{j=1}^\infty$ is a \emph{Cauchy sequence} in $M$ for every $t\in T$.
However, if $(M,d)$ is not complete, we may no longer infer the property
$\nu_n(f,g)\le\mu_n$, $n\in\NB$, of the pointwise limit~$f$ (as there may be no~$g$).

\smallbreak
(b) Condition \eq{e:limsup} is \emph{necessary\/} for the uniformly convergent
sequences $\{f_j\}$ and $\{g_j\}$: in fact, if $f_j\rra f$ and $g_j\rra g$ on $T$,
and $\nu_n(f,g)=o(n)$, then it follows from Lemma~\ref{l:ess}(d) that
  $\lim_{j\to\infty}\nu_n(f_j,g_j)=\nu_n(f,g)=o(n).$

\smallbreak
(c) Condition \eq{e:limsup} is `almost necessary' in the following sense. Suppose
$T\subset\RB$ is a measurable set with Lebesgue measure $\mathcal{L}(T)<\infty$,
$\{f_j\}$, $\{g_j\}\subset M^T$ are two sequences of measurable functions, which
converge pointwise (or almost everywhere) on $T$ to functions $f,g\in M^T$,
respectively, such that $\nu_n(f,g)=o(n)$. By Egorov's Theorem, given $\vep>0$,
there exists a measurable set $E_\vep\subset T$ such that
$\mathcal{L}(T\setminus E_\vep)\le\vep$, $f_j\rra f$ and $g_j\rra g$ on $E_\vep$.
So, as in the previous remark (b), we have
  $$\lim_{j\to\infty}\nu_n(f_j,g_j;E_\vep)=\nu_n(f,g;E_\vep)\le\nu_n(f,g)=o(n).$$
\end{rem}

\begin{exa} \label{ex:Di}
(a) Condition \eq{e:limsup} is \emph{not necessary\/} for the pointwise convergence
even if $g_j=c$ for all $j\in\NB$. To see this, let $T=[0,1]$ and $x,y\in M$ with
$x\ne y$. Given $j\in\NB$, define $f_j:T\to M$ by: $f_j(t)=x$ if $j!\cdot t$ is integer,
and $f_j(t)=y$ otherwise, $t\in[0,1]$. The pointwise precompact sequence
$\{f_j\}\subset M^T$ consists of bounded regulated functions (in fact,
$\nu_n(f_j,c)=o(n)$, and so, $f_j\in\mbox{\rm Reg}([0,1];M)=\mbox{\rm R}(c)$
for all $j\in\NB$). It converges pointwise on $T$ to the Dirichlet-type function
$\DC_{x,y}$ from \eq{e:Dxy}.~Note~that $\nu_n(\DC_{x,y},c)=nd(x,y)$, and so,
$\DC_{x,y}\notin\mbox{\rm R}(c)$. Since the usual Jordan variation $V(f_j)$ of $f_j$
on $T=[0,1]$ is equal to $2\cdot j!d(x,y)$, we find
  $$\nu_n(f_j,c)=d(x,y)\cdot\left\{
  \begin{array}{ccl}
\!n & \mbox{if} & n<2\cdot j!,\\[2pt]
\!2\cdot j! & \mbox{if} & n\ge 2\cdot j!,
  \end{array}\right.\quad\,\,n,j\in\NB.$$
Thus, $\lim_{j\to\infty}\nu_n(f_j,c)=d(x,y)\cdot n$, i.e., condition \eq{e:limsup}
does not hold.

\smallbreak
(b) Under the assumptions of Theorem~\ref{t:first}, condition \eq{e:limsup} does not
in general imply $\limsup_{j\to\infty}\nu_n(f_j,g)=o(n)$. To see this, let $g_j=f_j$ be
as in example (a) above, so that $g=\DC_{x,y}$. Given $n,j\in\NB$, choose a
collection $\{I_i\}_1^n\prec(0,1/j!)$ with $I_i=\{s_i,t_i\}$ such that $s_i$ is rational
and $t_i$ is irrational for all $i=1,\dots,n$. Noting that, by virtue of \eq{e:join},
  $$|(f_j,g)(\{s_i,t_i\})|=\sup_{z\in M}|d(y,z)-d(x,z)|=d(y,x),$$
we get
  $$\nu_n(f_j,g)\ge\sum_{i=1}^n|(f_j,g)(I_i)|=nd(y,x)\quad\mbox{for \,all}
     \quad n,j\in\NB.$$

\smallbreak
(c) The choice of an appropriate sequence $\{g_j\}$ is essential in Theorem~\ref{t:first}.
Let $\{x_j\}$, $\{y_j\}\subset M$ be two sequences, which converge in $M$ to
$x,y\in M$, respectively, $x\ne y$. Define $f_j:T=[0,1]\to M$ by
$f_j=\DC_{x_j,y_j}$, $j\in\NB$ (cf.\ \eq{e:Dxy}). Clearly, $\{f_j\}$ converges
uniformly on $T$ to $\DC_{x,y}$ (so, $\{f_j\}$ is pointwise precompact on $T$),
and $\nu_n(f_j,c)=nd(x_j,y_j)$ for all $n,j\in\NB$. Since
  $$|d(x_j,y_j)-d(x,y)|\le d(x_j,x)+d(y_j,y)\to0\quad\mbox{as}\quad j\to\infty,$$
we find $\lim_{j\to\infty}\nu_n(f_j,c)=nd(x,y)$, condition \eq{e:limsup} is not satisfied,
and Theorem~\ref{t:first} is inapplicable with $g_j=c$, $j\in\NB$.

\smallbreak
On the other hand, set $g_j=\DC_{x,y}$ for all $j\in\NB$. Given $\{s,t\}\subset T$,
we have, by virtue of \eq{e:fg2},
  $$|(f_j,g_j)(\{s,t\})|\le d(f_j(s),g_j(s))+d(f_j(t),g_j(t))\le2\vep_j,$$
where $\vep_j=\max\{d(x_j,x),d(y_j,y)\}\to0$ as $j\to\infty$. This implies the
inequality $\nu_n(f_j,g_j)\le2n\vep_j$, and so, condition \eq{e:limsup} is fulfilled.
\end{exa}

\section{Extensions of known selection theorems} \label{s:ext}

In this section, we consider extensions of two selection theorems from \cite{Fr} and
\cite{CyKe,Ko1}. The other selection theorems from the references in the Introduction
were shown to be particular cases of \cite{jmaa05}--\cite{MatTr}
(see Remark~\ref{r:5}).

\smallbreak
By virtue of the inequality $\nu_n(f_j,g_j)/n\le V(f_j,g_j)/n$, instead of condition 
\eq{e:limsup} in Theorem~\ref{t:first} we may assume that
$\limsup_{j\to\infty}V(f_j,g_j)<\infty$ or $\sup_{j\in\NB}V(f_j,g_j)<\infty$,
in which cases the resulting pointwise limit $f$ of a subsequence of $\{f_j\}$ satisfies the
regularity condition of the form $V(f,g)<\infty$.

\smallbreak
Making use of the notion of $\vep$-variation (Section~\ref{s:reg}), we get the following

\begin{thm} \label{t:eV}
Given $\varnothing\ne T\subset\RB$ and a metric space $(M,d)$, let
$\{f_j\}\subset M^T$ be a pointwise precompact sequence of functions such that
  \begin{equation} \label{e:eVf}
\limsup_{j\to\infty}V_\vep(f_j)<\infty\quad\mbox{for all}\quad\vep>0.
  \end{equation}
Then, there is a subsequence $\{f_{j_k}\}$ of $\{f_j\}$, which converges pointwise
on $T$ to a regulated function $f\in\mbox{\rm R}(c)$.
\end{thm}

\begin{pf}
Taking into account Theorem~\ref{t:first}, it suffices to verify that condition \eq{e:eVf}
implies $\limsup_{j\to\infty}\nu_n(f_j,c)\!=\!o(n)$, which is \eq{e:limsup} with $g_j\!=\!c$
for all~$j\!\in\!\NB$. In fact, by \eq{e:eVf}, for every $\vep>0$ there is
$j_0=j_0(\vep)\in\NB$ and $C(\vep)>0$ such that $V_\vep(f_j)\le C(\vep)$ for all
$j\ge j_0$. Definition \eq{e:evar} yields the existence
of \mbox{$g_j\in\mbox{\rm BV}(T;M)$} such that $d_\infty(f_j,g_j)\!\le\!\vep$
and $V(g_j)\le V_\vep(f_j)+(1/j)\le C(\vep)+1$ for all $j\ge j_0$. 
By \eq{e:Et} (where we replace $g_j$ and $g$ by $c$, and $f$---by $g_j$),
  $$\frac{\nu_n(f_j,c)}n\le2d_\infty(f_j,g_j)+\frac{\nu_n(g_j,c)}n
    \le2\vep+\frac{V(g_j)}n\le2\vep+\frac{C(\vep)+1}n$$
for all $j\ge j_0$ and $n\in\NB$. Consequently,
  $$\frac1n\limsup_{j\to\infty}\nu_n(f_j,c)\le\frac1n\sup_{j\ge j_0}\nu_n(f_j,c)
     \le2\vep+\frac{C(\vep)+1}n\quad\forall\,\vep>0,\,n\in\NB.$$
This implies that the left-hand side in this inequality tends to zero as \mbox{$n\to\infty$}:
given $\eta>0$, we set $\vep=\eta/4$ and choose a number $n_0=n_0(\eta)\in\NB$
such that $(C(\vep)+1)/n_0\le\eta/2$, which yields
$2\vep+(C(\vep)+1)/n\le\eta$ for all $n\ge n_0$.
\qed\end{pf}

\begin{rem} \label{r:4}
(a) If $M=\RB^N$ in Theorem~\ref{t:eV}, we may infer that $V_\vep(f)$ does not
exceed the limit superior from \eq{e:eVf}: in fact, it follows from
\cite[Proposition~3.6]{Fr} that $V_\vep(f)\le\liminf_{k\to\infty}V_\vep(f_{j_k})$
for all $\vep>0$.

\smallbreak
(b) Theorem~\ref{t:eV} extends Theorem~3.8 from \cite{Fr}, which has been
established for $T=[a,b]$ and $M=\RB^N$ under the assumption that
$\sup_{j\in\NB}V_\vep(f_j)<\infty$ for every $\vep>0$. The last assumption on
the uniform boundedness of $\vep$-variations is more restrictive than condition
\eq{e:eVf} as the following example shows.
\end{rem}

\begin{exa} \label{ex:fn}
Let $\{x_j\}$ and $\{y_j\}$ be two sequences from $M$ such that $x_j\ne y_j$ for all
$j\in M$ and, for some $x\in M$, $x_j\to x$ and $y_j\to x$ in $M$ as $j\to\infty$.
We set $f_j=\DC_{x_j,y_j}$, $j\in\NB$, and $f(t)=x$ for all $t\in T=[0,1]$.
The sequence $\{f_j\}\subset M^T$ converges uniformly on $T$ to the constant
function~$f$:
  $$d_\infty(f_j,f)=\max\{d(x_j,x),d(y_j,x)\}\to0\quad\mbox{as}\quad j\to\infty.$$
Given $\vep>0$, there is $j_0=j_0(\vep)\in\NB$ such that $d(x_j,y_j)\le\vep$
for all $j\ge j_0$, and so, by \eq{e:Vef}, $V_\vep(f_j)=0$ for all $j\ge j_0$, which
implies condition \eq{e:eVf}:
  $$\limsup_{j\to\infty}V_\vep(f_j)\le\sup_{j\ge j_0}V_\vep(f_j)=0.$$
On the other hand, if $k\in\NB$ is fixed and $0<\vep<d(x_k,y_k)/2$, then \eq{e:Vef}
gives $V_\vep(f_k)=\infty$, and so, $\sup_{j\in\NB}V_\vep(f_j)=\infty$.
\end{exa}

Now, we are going to present an extension of a Helly-type selection theorem from
\cite[Section~4, Theorem~1]{Ko1} and \cite[Theorem~2]{CyKe}.

\smallbreak
Let $\kappa:[0,1]\to[0,1]$ be a continuous, increasing and concave function such that
$\kappa(0)=0$, $\kappa(1)=1$, and $\kappa(\tau)/\tau\to\infty$ as $\tau\to+0$
(e.g., $\kappa(\tau)=\tau(1-\log\tau)$, $\kappa(\tau)=\tau^\alpha$ with $0<\alpha<1$,
or $\kappa(\tau)=1/(1-\frac12\log\tau)$, see \cite{Ko2}).

\smallbreak
Let $T=[a,b]$ be a closed interval in $\RB$, $a<b$. We set $|T|=b-a$, and if
$\{t_i\}_0^n\prec[a,b]$ is a partition of $T$ (i.e., $a=t_0<t_1<\dots<t_{n-1}<t_n=b$),
we also set $I_i=\{t_{i-1},t_i\}$ and $|I_i|=t_i-t_{i-1}$, $i=1,\dots,n$.

\smallbreak
The \emph{joint $\kappa$-variation\/} of functions $f,g\in M^T=M^{[a,b]}$
is defined by
  $$V_\kappa(f,g)\!=\!\sup\biggl\{
     \sum_{i=1}^n\!|(f,g)(I_i)|\biggl/\sum_{i=1}^n\!\kappa\bigl(|I_i|/|T|\bigr):
     \mbox{$n\!\in\!\NB$ and $\{t_i\}$}_0^n\!\prec\![a,b]\biggr\},$$
where $|(f,g)(I_i)|=|(f,g)(\{t_{i-1},t_i\})|$ is given by \eq{e:join}.

\smallbreak
Since $|(f,c)(I_i)|=d(f(t_{i-1}),f(t_i))$ is independent of a constant function
$c:[a,b]\to M$, the quantity $V_\kappa(f)\equiv V_\kappa(f,c)$ is the
\emph{Korenblum $\kappa$-variation\/} of $f\in M^{[a,b]}$, introduced in
\cite[p.~191]{Ko1} and \cite[Section~5]{Ko2} for $M=\RB$.

\smallbreak
The following theorem is a generalization of Theorem~2 from \cite{CyKe}, established
for real functions of bounded $\kappa$-variation under the assumption that
$\sup_{j\in\NB}V_\kappa(f_j)<\infty$ and based on the decomposition of any
$f\in\RB^{[a,b]}$ with $V_\kappa(f)<\infty$ into the difference of two real
$\kappa$-decreasing functions.

\begin{thm} \label{t:Ko}
Under the assumptions of Theorem~\ref{t:first}, suppose that condition\/
\eq{e:limsup} is replaced by the following one\/{\rm:}
$\limsup_{j\to\infty}V_\kappa(f_j,g_j)<\infty$.
Then, there is a subsequence of $\{f_j\}$, which converges pointwise on $T=[a,b]$
to a function $f\in\mbox{\rm R}(g)$ such that $V_\kappa(f,g)<\infty$.
\end{thm}

\begin{pf}
In order to show that \eq{e:limsup} is satisfied, let $n\in\NB$,
$\{I_i\}_1^n\prec[a,b]$ with $I_i=\{s_i,t_i\}$, and set $I_i'=\{t_i,s_{i+1}\}$ and
$|I_i'|=s_{i+1}-t_i$, $i=1,\dots,n-1$. By the definition of $V_\kappa(f_j,g_j)$ and
the concavity of function $\kappa$, we have
  \begin{align*}
\sum_{i=1}^n\!|(f_j,g_j)(I_i)|&\le\!|(f_j,g_j)(\{a,s_1\})|+\sum_{i=1}^n\!|(f_j,g_j)(I_i)|
  \\[-4pt]
&\qquad+\sum_{i=1}^{n-1}\!|(f_j,g_j)(I_i')|+|(f_j,g_j)(\{t_n,b\})|\\[-2pt]
&\le\!\biggl[\!\kappa\biggl(\!\frac{s_1\!-\!a}{|T|}\!\biggr)
  \!+\!\sum_{i=1}^n\!\!\kappa\biggl(\frac{|I_i|}{|T|}\biggr)
  \!+\!\sum_{i=1}^{n-1}\!\!\kappa\biggl(\frac{|I_i'|}{|T|}\biggr)
   \!+\!\kappa\biggl(\!\frac{b\!-\!t_n}{|T|}\!\biggr)\!\biggr]V_\kappa(f_j,g_j)\\[-2pt]
&\le(2n\!+\!1)\kappa\biggl(\frac{1}{(2n\!+\!1)(b-a)}\biggl[(s_1\!-\!a)
  +\sum_{i=1}^n(t_i\!-\!s_i)\\[-4pt]
&\qquad+\sum_{i=1}^{n-1}(s_{i+1}\!-\!t_i)
  +(b-t_n)\biggr]\biggr)V_\kappa(f_j,g_j)\\[-4pt]
&\le(2n\!+\!1)\kappa\biggl(\frac1{2n\!+\!1}\biggr)V_\kappa(f_j,g_j).
  \end{align*}
Thus,
  $$\frac{\nu_n(f_j,g_j)}n\le\biggl(2+\frac1n\biggr)\kappa\biggl(\frac1{2n+1}\biggr)
     V_\kappa(f_j,g_j)\quad\mbox{for all}\quad j,n\in\NB,$$
and so, condition \eq{e:limsup} in Theorem~\ref{t:first} is satisfied.

\smallbreak
Let $f\in\mbox{\rm R}(g)$ be the pointwise limit of a subsequence $\{f_{j_m}\}$
of $\{f_j\}$. Arguing as in the proof of Lemma~\ref{l:ess}(c), we get
  $$V_\kappa(f,g)\le\liminf_{m\to\infty}V_\kappa(f_{j_m},g_{j_m})
     \le\limsup_{j\to\infty}V_\kappa(f_j,g_j)<\infty.\qquad\Box$$
\end{pf}

\begin{rem} \label{r:5}
Since Theorem~\ref{t:first} is an extension of results from \cite{jmaa05,Ischia},
it also contains as particular cases all pointwise selection principles based on various
notions of generalized variation. These principles may be further
generalized in the spirit of Theorem~\ref{t:Ko} replacing the increment
$|f(I_i)|=d(f(s_i),f(t_i))$ applied in \cite{jmaa05}--\cite{MatTr} by the joint
increment $|(f,g)(I_i)|$ from \eq{e:join}.
\end{rem}

Finally, it is worth mentioning that the joint modulus of variation \eq{e:jmv}, defined
by means of \eq{e:join}, plays an important role in the extension of a result from
\cite{jmaa08} to metric space valued functions:

\begin{thm} \label{t:ChMa}
Given $\varnothing\ne T\subset\RB$ and a metric space $(M,d)$, let 
$\{f_j\}\subset M^T$ be a pointwise precompact sequence of functions such that
  $$\lim_{N\to\infty}\sup_{j,k\ge N}\nu_n(f_j,f_k)=o(n).$$
Then, there is a subsequence of $\{f_j\}$, which converges pointwise on~$T$.
\end{thm}

Taking into account Lemmas \ref{l:pjin} and \ref{l:ess}, the proof of this theorem
follows the same lines as the ones given in \cite[Theorem~1]{jmaa08} (where
$(M,d,+)$ is a metric semigroup and $\nu_n(f_j,f_k)$ is defined via \eq{e:mes}),
and so, it is omitted. Note that the limit of a pointwise convergent subsequence
of $\{f_j\}$ in Theorem~\ref{t:ChMa} may be a non-regulated function.
For more details, examples and relations with previously known `regular' and
`irregular' versions of pointwise selection principles  from \cite{Piazza, Schrader, MZ08}
we refer to \cite[Section~5.2]{MMS}, \cite{jmaa08,waterman80,Manisc}.

\section{Conclusions}

In the context of functions of one real variable taking values in a metric space,
we have presented a pointwise selection principle, which can be applied to arbitrary
(regulated and non-regulated) sequences of functions. It is based on notions of
joint increment and joint modulus of variation, the latter being a nondecreasing
sequence of pseudometrics on the appropriate product space. Our selection principle
extends the classical Helly Theorem and contains as particular cases many selection
theorems based on various notions of generalized variation of functions.
In contrast to previously established selection principles, the main assumption
in our theorem on the limit superior is `almost necessary', and it is shown by
examples that in a certain sense this assumption is sharp. The notion of
conditionally regulated functions explains the limitations of previously known
selection results.

\end{document}